\documentclass[11pt,a4paper]{article}
\usepackage{amsfonts}
\usepackage{}
\usepackage[pagewise]{lineno}
\usepackage{graphicx}
\usepackage{amsfonts,amsmath, amssymb}

\newtheorem{theo}{\textbf{\ \ \quad Theorem}}[section]

\newtheorem{remark}{\textbf{\ \ \quad Remark}}[section]
\newtheorem{col}{\textbf{\ \ \quad Corollary}}[section]
\newtheorem{prop}{\textbf{\ \ \quad Proposition}}[section]

\newcommand{\lbl}[1]{\label{#1}}

\newcommand{\be}{\begin{equation}}
\newcommand{\ee}{\end{equation}}
\newcommand\bes{\begin{eqnarray}}
\newcommand\ees{\end{eqnarray}}
\newcommand{\bess}{\begin{eqnarray*}}
\newcommand{\eess}{\end{eqnarray*}}

{\setlength\arraycolsep{2pt}}

\title{Blowup solutions for stochastic parabolic equations }

\author{Guangying Lv$^{a,b}$,  Jinlong Wei$^c$\\
\\
\ \\
   {\small \it $^a$ Institute of Applied Mathematics, Henan University, Kaifeng, Henan 475001, China}\\
  {\small \it $^b$ Center for Applied Mathematics, Tianjin University, Tianjin 300072, China }\\
   {\small \tt gylvmaths@henu.edu.cn}\\
   {\small \it $^c$ School of Statistics and Mathematics, Zhongnan University of}\\
   {\small \it Economics and Law, Wuhan 430073, China}\\
    {\small \tt  weijinlong.hust@gmail.com }
}

\begin{document}
\maketitle

\medskip

\begin{abstract}
In this short paper, we are concerned with the  blowup phenomenon of stochastic parabolic equations.
By using comparison principle and the results of deterministic
parabolic equations, we obtain blowup results of solutions for stochastic parabolic
equations.

{\bf Keywords}: Blowup; Stochastic parabolic equation; Comparison principle.

AMS subject classifications (2010): 35K20, 60H15, 60H40.

\end{abstract}

\baselineskip=15pt

\section{Introduction}
\setcounter{equation}{0}

For deterministic partial differential equations,
finite time blowup phenomenon has been studied by many authors, see the book \cite{Hubook2018}.
There are two cases to study this problem. One is bounded domain and the other is whole space.
On the bounded domain, the $L^p$-norm of solutions ($p>1$) will blow up in finite time.
The other is  the whole space. The following 'Fujita Phenomenon' has been attracted
in the literature. Consider the following Cauchy problem
   \bes\left\{\begin{array}{lll}
u_t=\Delta u+u^p, \quad t>0,\ \ p>0, &x\in\mathbb{R}^d,\\
u(0,x)=u_0(x), \ \ \ &x\in\mathbb{R}^d.
   \end{array}\right.\lbl{1.1}\ees
It has been proved that:
 \begin{quote}
(i) if $0<p<1$, then every nonnegative solution is global, but not necessarily
unique;

(ii) if $1<p\leq1+\frac{2}{d}$, then any nontrivial, nonnegative solution
blows up in finite time;

(iii) if $p>1+\frac{2}{d}$, then $u_0\in\mathcal{U}$ implies that
$u(t,x,u_0)$ exists globally;

(iv) if $p>1+\frac{2}{d}$, then $u_0\in\mathcal{U_1}$ implies that
$u(t,x,u_0)$ blows up in finite time,
 \end{quote}
where $\mathcal{U}$ and $\mathcal{U_1}$ are defined as follows
   \bess
&\mathcal{U}=\left\{v(x)|v(x)\in BC(\mathbb{R}^d,\mathbb{R}_+),
v(x)\leq \delta e^{-k|x|^2},\ k>0,\delta=\delta(k)>0\right\},\\
&\mathcal{U_1}=\left\{v(x)|v(x)\in BC(\mathbb{R}^d,\mathbb{R}_+),
v(x)\geq c e^{-k|x|^2},\ k>0,c\gg1\right\}.
  \eess
Here $BC=\{$ bounded and uniformly continuous functions $\}$, see
Fujita \cite{F1966,F1970}.

It is easy to see that for the whole space, there are four types of behaviors for problem (\ref{1.1}), namely,
(1) global existence unconditionally but uniqueness fails in certain solutions, (2) global existence with restricted
initial data, (3) blowing up unconditionally, and (4) blowing up with restricted initial data. The
occurrence of these behaviors depends on the combination effect of the nonlinearity represented by the parameter $p$, the size of the initial datum $u_0(x)$, represented by the choice of $\mathcal{U}$
or $\mathcal{U_1}$, and the dimension of the space.

Now, we recall some known results of stochastic partial differential equations (SPDEs).
In this paper, we only focus on the stochastic parabolic equations.
It is known that the existence and uniqueness of
global solutions to SPDEs can  be   established under appropriate conditions (\cite{DW2014}).
For the finite time blowup phenomenon of stochastic parabolic equations, we first consider
the case on a bounded domain. Consider the following equation
   \bes\left\{\begin{array}{lll}
   du=(\Delta u+f(u))dt+\sigma(u)dW_t, \ \ \qquad t>0,&x\in  D,\\[1.5mm]
   u(x,0)=u_0(x)\geq0, \ \ \ &x\in D,\\[1.5mm]
   u(t,x)=0, \qquad \qquad \qquad \qquad \qquad \qquad t>0,  &x\in\partial D.
    \end{array}\right.\lbl{1.2}\ees
Da Prato-Zabczyk \cite{PZ1992} considered  the existence of global solutions of (\ref{1.2})  with
additive noise ($\sigma$ is a constant).
Dozzi and L\'{o}pez-Mimbela
\cite{DL2010} studied   equation (\ref{1.2}) with
$\sigma(u)=u$ and  proved that if $f(u)\geq u^{1+\alpha}$ ($\alpha>0$) and the initial data is large enough, the solution will blow up in finite time,
and that if $f(u)\leq u^{1+\beta}$ ($\beta$ is a certain positive constant) and the initial data is
small enough, the solution will exist globally. A natural question arises: If
$\sigma$ does not satisfy the global Lipschitz condition, what can we say about the solution?
Will it blow up in finite time or exist globally? Chow \cite{C2011} answered part of this question.
Lv-Duan \cite{LD2015} described the competition between the nonlinear term and noise term for equation (\ref{1.2}).
Bao-Yuan \cite{BY2014} and Li et al.\cite{LPJ2016} obtained the existence of local solutions of (\ref{1.2}) with jump process
and L\'{e}vy process, respectively. For blowup phenomenon of stochastic functional
parabolic equations, we refer to \cite{CL2012,FP2015}.
 In a somewhat different case, Mueller \cite{M1991} and, later,
Mueller-Sowers \cite{MuS1993} investigated the problem of a noise-induced explosion for a
special case of equation (\ref{1.2}), where $f(u)\equiv0,\,\sigma(u)=u^\gamma$ with
$\gamma>0$ and $W(t,x)$ is a space-time white noise. It was shown that the solution will
explode in finite time with positive probability for some $\gamma>3/2$.

For the whole space, Foondun et al. \cite{FLN2018} considered the nonexistence of global solutions for the Cauchy problem of stochastic fractional parabolic equations.
Comparing with the deterministic parabolic equations, they only obtained
the result similar to type (4), also see \cite{Wang2019}.
In paper \cite{LW2019}, we established the similar results to types (1)
and (3). The method used there is the properties of heat kernel.
But in \cite{LW2019}, we only obtained the results for one dimension in the
whole space. More precisely, we obtained the following result.
   \begin{prop}\lbl{p1.1}
Suppose $\sigma^2(u,x,t)\geq C_0u^{4}$, $C_0>0$, then the solutions of the following
equation will blow up in finite
time for any nontrivial nonnegative initial data $u_0$,
  \bes\left\{\begin{array}{llll}
du=\Delta udt+\sigma(u,x,t)dB_t,\ \ t>0,\ &x\in\mathbb{R},\\
u(x,0)=u_0(x)\gneqq0, \ \  &x\in\mathbb{R}.
   \end{array}\right. \lbl{1.3}\ees
  \end{prop}

In this paper, we consider the blowup phenomenon of SPDEs included fractional diffusion
equations, and generalize of the result Proposition \ref{p1.1}. There are a lot of work to do for
SPDEs comparing with the deterministic case, see \cite{Yang2019}.

This paper is arranged as follows. In Sections 2,  we state out the main results and the proofs. Throughout this paper, we write $C$ as a general positive constant and $C_i$, $i=1,2,\cdots$ as
a concrete positive constant.
\section{Main results and proofs}
\setcounter{equation}{0}
In this section, we recall some known results and state out the main results.
\begin{prop}\lbl{p2.1}\cite{Pinsky1997}
Consider the following equation
\bes\left\{\begin{array}{lll}
u_t=\Delta u+a(x)u^p, \ \ t\in(0,T), &x\in  \mathbb{R}^d,\\[1.5mm]
u(0,x)=u_0(x)\geq0, \ \ \ &x\in \mathbb{R}^d,
\end{array}\right.\lbl{2.1}\ees
where $0\lneqq a\in C^\beta(\mathbb{R}^d)$, $\beta>0$ and $p>1$.

1. Assume that for large $|x|$ and constants
$c_1,c_2>0$, $c_1|x|^m\leq a(x)\leq c_2|x|^m$, $m>-1$. If $1<p\leq1+(2+m)/d$, then (\ref{2.1}) does not possess global solutions, for any choice of initial data $u_0\gneqq0$;

2. Assume that for large $|x|$ and constants
$c_2>0$, $a(x)\leq c_2|x|^{-2}$. If $d=1$ and $1<p\leq2$, then (\ref{2.1}) does not possess global solutions, for
any choice of initial data $u_0\gneqq0$.
\end{prop}

We now consider the following Cauchy problem
   \bes \left\{\begin{array}{llll}
du=\Delta udt+b(t)|u|^pdW(t,x),\ \ t>0,\ &x\in\mathbb{R},\\
u(0,x)=u_0(x)\gneqq0, \ \  &x\in\mathbb{R},
   \end{array}\right.\lbl{2.2}\ees
where $W(t,x)$ is a white noise both in time and space, and $0< b\in C(\mathbb{R}_+)$.
In the rest of paper, we always assume that the initial data is a nonnegative continuous function.
A mild solution to (\ref{2.2}) in sense of Walsh \cite{walsh1986} is any
$u$ which is adapted to the filtration generated by the
white noise and satisfies the following evolution equation
   \bess
u(t,x)=\int_{\mathbb{R}}K(t,x-y)u_0(y)dy
+\int_0^t\int_{\mathbb{R}}K(t-s,x-y)b(s)|u|^p(s,y)W(ds,dy),
   \eess
where $K(t,x)$ denotes the heat kernel of Laplacian operator.
We get the following results.
   \begin{theo}\lbl{t2.1} Let $0< b\in C(\mathbb{R}_+)$ such that
\bess
b_\infty:=\int_0^\infty (b^{-1}(t))^2dt<+\infty.
   \eess
If $1<p\leq2$, then (\ref{2.2}) does not possess global solutions, for
any choice of initial data $u_0\gneqq0$.
\end{theo}

{\bf Proof.}
By taking the second moment and using the Walsh isometry,
we get for any $0<t\leq T$ ($T>0$ is any fixed number)
   \bess
&&\mathbb{E}|u(t,x)|^2
\\&=&\left(\int_{\mathbb{R}}K(t,x-y)u_0(y)dy\right)^2
+\int_0^tb^2(s)\int_{\mathbb{R}}K^2(t-s,x-y) \mathbb{E}|u|^{2p}(s,y)dyds\\
&\geq&\left(\int_{\mathbb{R}}K(t,x-y)u_0(y)dy\right)^2
+\int_0^tb^2(s)\int_{\mathbb{R}}K^2(t-s,x-y) \left[\mathbb{E}|u|^{2}(s,y)\right]^pdyds.
   \eess

Let $0\lneqq a\in C^\beta(\mathbb{R})$ such that
for large $|x|$, there is a positive constant $c_1>0$ satisfying $a(x)\leq c_1|x|^{-2}$. Then $a_\infty:=\int_{\mathbb{R}}a^2(y)dy<+\infty$. Observing that
 \bess
&& \left(\int_0^t\int_{\mathbb{R}}K(t-s,x-y) a(y)\left[\mathbb{E}|u|^2(s,y)\right]^{\frac{p}{2}}dyds\right)^2
\\ &\leq&
 \left(\int_0^t\left[\int_{\mathbb{R}}K^2(t-s,x-y) \left[\mathbb{E}|u|^2(s,y)\right]^{p}dy\right]^{\frac12}
 \left[\int_{\mathbb{R}}a(y)^2dy\right]^{\frac12}ds\right)^2
 \\ &\leq&a_\infty \int_0^tb^2(s)\int_{\mathbb{R}}K^2(t-s,x-y) \left[\mathbb{E}|u|^2(s,y)\right]^{p}dyds\int_0^t\frac{1}{b^2(s)}ds
\\ &\leq&a_\infty b_\infty \int_0^tb^2(s)\int_{\mathbb{R}}K^2(t-s,x-y) \left[\mathbb{E}|u|^2(s,y)\right]^{p}dyds.
\eess

We deduce that
\bess
&&\mathbb{E}|u(t,x)|^2
\\&\geq&\left(\int_{\mathbb{R}}K(t,x-y)u_0(y)dy\right)^2
+\frac{1}{a_\infty b_\infty }\left(\int_0^t\int_{\mathbb{R}}K(t-s,x-y) a(y)\left(\mathbb{E}|u|^2(s,y)\right)^{\frac{p}{2}}dyds\right)^2\\
&\geq&\frac{1}{2(a_\infty b_\infty \vee1)}\left(\int_{\mathbb{R}}K(t,x-y)u_0(y)dy\right. \left.+\int_0^t\int_{\mathbb{R}}K(t-s,x-y) a(y)\left(\mathbb{E}|u|^2(s,y)\right)^{\frac{p}{2}}dyds\right)^2.
\eess

For $0<t\leq T$, we set $v(t,x)=\left(\mathbb{E}|u(t,x)|^2\right)^{\frac{1}{2}}$,
then $v$ is an upper solution of the following equation
\bes\left\{\begin{array}{lll}
v_t=\Delta v+\frac{1}{\sqrt{2(a_\infty b_\infty \vee1)}}a(x)v^p, \ \ \ t\in(0,T),&x\in  \mathbb{R},\\[1.5mm]
v(0,x)=\frac{1}{\sqrt{2(a_\infty b_\infty \vee1)}}u_0(x)\geq0, \ \ \ &x\in \mathbb{R}.
\end{array}\right.\lbl{2.3}\ees
By Proposition \ref{p2.1}, if $1<p\leq 2$, then (\ref{2.3}) does not possess global solutions, for
any choice of initial data $u_0\gneqq0$. So does the equation (\ref{2.2}).
The proof is complete. $\Box$

\begin{remark}
In deterministic parabolic equations, we will assume that $d\geq1$, which is
different from the stochastic case. The main reason is that the It\^{o} isometry
may be not right for $d>1$.
  \end{remark}

Now, we extend Theorem \ref{t2.1} to stochastic fractional parabolic equations. Initially, let us consider  the following deterministic equation:
\bes\left\{\begin{array}{lll}
u_t=(-\Delta)^{\frac{\alpha}{2}} u+u^p, \ \ \qquad t>0,&x\in  \mathbb{R}^d,\\[1.5mm]
   u(0,x)=u_0(x)\geq0, \ \ \ &x\in \mathbb{R}^d.
    \end{array}\right.\lbl{2.4}\ees
Sugitani \cite{Sug1975} proved that if
$p\in(1,1+\alpha/d]$ and the initial data is non-trivial and non-negative, the solution of (\ref{2.4}) will blow up in finite time. Next, let us consider the following stochastic fractional Laplacian equation
  \bes \left\{\begin{array}{llll}
du=-(-\Delta)^{\frac\alpha2} udt+b(t)a(x)|u|^pdW(t,x),\ \ t>0,\ &x\in\mathbb{R},\\
u(0,x)=u_0(x)\gneqq0, \ \  &x\in\mathbb{R},
   \end{array}\right.\lbl{2.5}\ees
where $W(t,x)$ is described in (\ref{2.1}). We will study the nonexistence of global solutions to (\ref{2.5}). Here, a solution to (\ref{2.5}) is defined by
   \bess
u(t,x)=\int_{\mathbb{R}}p(t,x-y)u_0(y)dy
+\int_0^t\int_{\mathbb{R}}b(s)a(y)p(t-s,x-y)|u|^p(s,y)W(ds,dy),
   \eess
where $p(t,x)$ denotes the heat kernel of fractional Laplacian operator. By a similar discussion as in Theorem \ref{t2.1}, we have

 \begin{col}\lbl{c2.1}
Assume $0<\alpha<2$ and $1<p\leq1+\alpha$. Assume further that
   \bess
b_\infty:=\int_0^\infty (b^{-1}(t))^2dt<+\infty, \ \ a_\infty:=\int_0^\infty (a^{-1}(x))^2dx<+\infty.
   \eess
Then (\ref{2.5}) does not possess global solutions, for
any choice of initial data $u_0\gneqq0$.
\end{col}

Next, we consider
the following stochastic parabolic equation
  \bes \left\{\begin{array}{llll}
du_t=\Delta udt+b(t)a(x)|u|^pdB_t,\ \ t>0,\ &x\in\mathbb{R}^d,\\
u(x,0)=u_0(x)\gneqq0, \ \  &x\in\mathbb{R}^d,
   \end{array}\right.\lbl{2.6}\ees
where $B_t$ is one-dimensional Brownian motion, $a\in C(\mathbb{R}^d)$ and
$0< b\in C(\mathbb{R}_+)$.
A mild solution to (\ref{2.6}) in sense of Walsh \cite{walsh1986} is any
$u$ which is adapted to the filtration generated by the
white noise and satisfies the following evolution equation
   \bess
u(t,x)=\int_{\mathbb{R}^d}K(t,x-y)u_0(y)dy
+\int_0^t\int_{\mathbb{R}^d}K(t-s,x-y)a(y)b(s)|u|^p(s,y)dydB_s,
   \eess
where $K(t,x)$ denotes the heat kernel of Laplacian operator.

\begin{theo}\lbl{t2.2} Let $p\geq2$ and for large $|x|$, $c_1|x|^m\leq |a(x)|\leq c_2|x|^m$ with $c_1,c_2>0,m\in\mathbb{R}$. Assume that $m>d-2$ and $2\leq p\leq 1+(2+m)/d$, that
   \bess
b_\infty:=\int_0^\infty (b^{-1}(t))^2dt<+\infty.
   \eess
Then (\ref{2.6})  does not possess global solutions for any choice of initial data $u_0\gneqq0$.
  \end{theo}

{\bf Proof.} By taking the second moment and using the Walsh isometry,
we get for any $0<t\leq T$ ($T>0$ is any fixed number)
   \bess
&&\mathbb{E}|u(t,x)|^2
\\&=&\left(\int_{\mathbb{R}^d}K(t,x-y)u_0(y)dy\right)^2
+\int_0^tb^2(s)\left(\int_{\mathbb{R}^d}K(t-s,x-y) a(y)\mathbb{E}|u|^p(s,y)dy\right)^2ds\\
&\geq&\left(\int_{\mathbb{R}^d}K(t,x-y)u_0(y)dy\right)^2
+\frac{1}{b_\infty}\left(\int_0^t\int_{\mathbb{R}^d}K(t-s,x-y) a(y)\left(\mathbb{E}|u|^2(s,y)\right)^{\frac{p}{2}}dyds\right)^2\\
&\geq&\frac{1}{2(b_\infty\vee1)}\left(\int_{\mathbb{R}^d}K(t,x-y)u_0(y)dy\right.
\left.+\int_0^t\int_{\mathbb{R}^d}K(t-s,x-y) a(y)\left(\mathbb{E}|u|^2(s,y)\right)^{\frac{p}{2}}dyds\right)^2,
   \eess
where we used the fact that
   \bess
\int_0^tf(s)ds\leq\left(\int_0^t(b^{-1}(s))^2ds\right)^{\frac{1}{2}}\left(\int_0^t(b(s)f(s))^2ds\right)^{\frac{1}{2}}.
   \eess
Set
   \bess
v(t,x)=\left(\mathbb{E}|u(t,x)|^2\right)^{\frac{1}{2}},\ \ \
0<t\leq T,
   \eess
then $v(t,x)$ is an upper solution of the following equation
\bes\left\{\begin{array}{lll}
v_t=\Delta u+\frac{1}{\sqrt{2(b_\infty\vee1)}}a(x)v^p, \ \ \qquad t>0,&x\in  \mathbb{R}^d, \ t\in(0,T),\\[1.5mm]
v(0,x)=\frac{1}{\sqrt{2(b_\infty\vee1)}}u_0(x)\geq0, \ \ \ &x\in \mathbb{R}^d.
\end{array}\right.\lbl{2.7}\ees
Under the condition of Proposition \ref{p2.1} that the solution of (\ref{2.7})  does not possess global solutions for any choice of initial data $u_0\gneqq0$,
so does the equation (\ref{2.6}). $\Box$

\begin{remark}\lbl{r2.1}

(i) It follows from Theorem \ref{t2.2} that we consider the equation (\ref{2.7}) in any dimension, which
is different from that of \cite{LW2019}.

(ii) About the stochastic Fujita index, we have the following example. Consider the following equation
  \bes\left\{\begin{array}{llll}
du=\Delta udt+(1+t)^qa(x)|u|^pdB_t,\ \ t>0,\ &x\in\mathbb{R}^d,\\
u(x,0)=u_0(x)\gneqq0, \ \  &x\in\mathbb{R}^d.
   \end{array}\right. \lbl{2.8}\ees
Theorem \ref{t2.2} shows that if $2\leq p\leq1+(m+2)/d$, then (\ref{2.8}) does not possess global solutions for
any choice of initial data $u_0\gneqq0$.
\end{remark}

Lastly, we consider the impact of additive noise on parabolics. That is to say, we consider the following
Cauchy problem
  \bes \left\{\begin{array}{llll}
du=\Delta udt+|u|^pdt+\sigma(t,x)dB_t,\ \ t>0,\ &x\in\mathbb{R}^d,\\
u(x,0)=u_0(x)\gneqq0, \ \  &x\in\mathbb{R}^d.
   \end{array}\right.\lbl{2.9}\ees
We have the following result.
\begin{theo}\lbl{t2.3} (i) If $d=1$ and $2\leq p\leq 3$ and the initial data satisfies
   \bes
\left(\int_{\mathbb{R}^d}K(t,x-y)u_0(y)dy\right)^2
\geq\frac{4}{1-\beta}\int_0^t\left(\int_{\mathbb{R}^d}
K(t-s,x-y)\sigma(s,y)dy\right)^2ds, \ \
\forall \ t>0, \lbl{2.10}
  \ees
for some real number $\beta\in (0,1)$. Then the solution of (\ref{2.9}) will blow up in finite time.

(ii) For general $d$, we assume that $u_0\in \mathcal{U_1}$ and (\ref{2.10}) holds. If $p>1+d/2$ and $p\geq 2$, then the solution of (\ref{2.9}) will blow up in finite time.
  \end{theo}

{\bf Proof.} By taking the second moment,
we get for any $0<t\leq T$ ($T>0$ is any fixed number)
   \bes
&&\mathbb{E}|u(t,x)|^2
\nonumber\\ &=&\mathbb{E}\left(\int_{\mathbb{R}^d}K(t,x-y)u_0(y)dy
+\int_0^t\int_{\mathbb{R}^d}K(t-s,x-y)u^p(s,y)dyds\right.
\nonumber\\
&&\left.
+\int_0^t\int_{\mathbb{R}^d}K(t-s,x-y) \sigma(s,y)dydB_s\right)^2\nonumber\\
&\geq&\beta\left(\int_{\mathbb{R}^d}K(t,x-y)u_0(y)dy+\int_0^t\int_{\mathbb{R}^d}K(t-s,x-y)
[\mathbb{E}|u|^2]^{p/2}dyds\right)^2\nonumber\\
&&
-\frac{\beta}{1-\beta}\int_0^t\left(\int_{\mathbb{R}^d}K(t-s,x-y) \sigma(s,y)dy\right)^2ds\nonumber\\
&\geq&\beta\left(\int_{\mathbb{R}^d}K(t,x-y)u_0(y)dy-\frac{1}{\sqrt{1-\beta}}\left[\int_0^t\left(\int_{\mathbb{R}^d}K(t-s,x-y) \sigma(s,y)dy\right)^2ds\right]^{1/2}\right.\nonumber\\
&&\left.+\int_0^t\int_{\mathbb{R}^d}K(t-s,x-y)\left(\mathbb{E}|u|^2(s,y)\right)^{\frac{p}{2}}dyds\right)^2,
   \lbl{2.11}\ees
where we used
   \bess
(a+b)^2\geq \beta a^2-\frac{\beta}{1-\beta}b^2, \ \ \ a,\ b\in\mathbb{R}.
   \eess
In view of (\ref{2.10}), from (\ref{2.11}), we arrive at
\bes
\mathbb{E}|u(t,x)|^2
\geq \frac{\beta}{4}\left(\int_{\mathbb{R}^d}K(t,x-y)u_0(y)dy
+\int_0^t\int_{\mathbb{R}^d}K(t-s,x-y)
\left(\mathbb{E}|u|^2(s,y)\right)^{\frac{p}{2}}dyds\right)^2.
\lbl{2.12}\ees

From (\ref{2.12}), by using the comparison principle, we get the desired result.  $\Box$

\begin{remark}\lbl{r2.3} It follows from Theorem \ref{t2.3} that
the additive noise prevents the finite time blowup.

(ii) One can use a similar method to deal with the case of bounded domain if the operator admits a heat kernel. For references in this direction for the deterministic equations, we refer to \cite{M1990}.
  \end{remark}

\medskip

\noindent {\bf Acknowledgment} The first author was supported in part by NSFC of China grants 11771123, 11501577.


\begin{thebibliography}{99}\label{ref:ref}
\addtolength{\itemsep}{-1.4ex}
\bibitem{BY2014} J. Bao and C. Yuan, {\em Blow-up for stochastic reactin-diffusion equations with jumps},
J Theor. Probab. {\bf29} (2016) 617-631.





\bibitem{C2011} P-L. Chow, {\em  Explosive solutions of stochastic reaction-diffusion
equations in mean $L^p$-norm}, J. Differential Equations {\bf 250} (2011) 2567-2580.


\bibitem{CL2012} P-L. Chow and K. Liu, {\em  Positivity and explosion in mean $L^p$-norm
of stochastic functional parabolic equations of retarded type}, Stochastic Process. Appl., {\bf 122} (2012) 1709-1729.


\bibitem{DW2014} J. Duan and W. Wang, {\em Effective Dynamics of Stochastic
Partial Differential Equations}, Elsevier, 2014.

\bibitem{DL2010} M. Dozzi and J. A. L\'{o}pez-Mimbela, {\em
Finite-time blowup and existence of global positive solutions of a semi-linear spde},
Stochastic Process. Appl., {\bf120} (2010) 767-776.



\bibitem{FP2015} M. Foondun and Rana D. Parshad, {\em On non-existence of global solutions to a
class of stochastic heat equations}, Proc. Amer. Math. Soc. {\bf143} (2015) 4085-4094.

\bibitem{FLN2018} M. Foondun, W. Liu and E. Nane, {\em Some non-existence
results for a class of stochastic partial differential equations}, J. Differential Equations {\bf266} (2019) 2575¨C2596.

\bibitem{F1966} H. Fujita, {\em On the blowing up of solutions of the Cauchy problem for
$u_t-\Delta u=u^{1+\alpha}$}, J. Fac. Sci. Univ. Tokyo Sect.
IA Math. {\bf 13} (1966) 109-124.

\bibitem{F1970} H. Fujita, {\em On some nonexistence and nonuniqueness theorems for nonlinear parabolic equations},
 Proc. Symp. Pure Math. {\bf XVIII} (1970) 105-113.


\bibitem{Hubook2018} B. Hu, {\em Blow-up Theories for Semilinear Parabolic Equations},
Lecture Notes in Mathematics ISSN print edition: 0075-8434, Springer Heidelberg Dordrecht London New York, 2018.




\bibitem{LPJ2016} K. Li, J. Peng and J. Jia, {\em Explosive solutions
of parabolic stochastic partial differential equations with L$\acute{e}$vy
noise}, arXiv:1306.01676.





\bibitem{LD2015} G. Lv and J.   Duan, {\em Impacts of Noise on a Class of Partial Differential Equations},
J. Differential Equations, {\bf258} (2015) 2196-2220.


\bibitem{LW2019} G. Lv and J. Wei, {\em Global existence and non-existence of stochastic parabolic equations}, arXiv:1902.07389.



\bibitem{M1990} P. Meier, {\em On the critical exponent for reaction-diffusion equations},
Arch. Rational Mech. Anal.  {\bf109}  (1990) 63-71.


\bibitem{M1991} C. Mueller, {\em Long time existence for the heat equation with a noise term},
Probab. Theory Related Fields {\bf 90} (1991) 505-517.

\bibitem{MuS1993} C. Mueller and R. Sowers, {\em Blowup for the heat equation with a noise term},
Probab. Theory Related Fields {\bf 93} (1993) 287-320.



\bibitem{Pinsky1997} R. Pinsky, {\em Existence and Nonexistence of global
solutions for $u_t=\Delta u+a(x)u^p$ in $\mathbb{R}^d$},
J. Differential Equations {\bf 133} (1997) 152-177.


\bibitem{PZ1992} G. Da Prato and J. Zabczyk, {\em Nonexplosion, boundedness and ergodicity
for stochastic semilinear equations}, J. Differential Equations {\bf 98} (1992) 181-195.





\bibitem{Sug1975} S. Sugitani, {\em On nonexistence of global
solutions for some nonlinear integral equations}, Osaka J. Math.,
{\bf12} (1975) 45-51.




\bibitem{walsh1986} J.B. Walsh, {\em
An introduction to Stochastic Partial Differential Equations},
volume 1180 of Lecture Notes in Math., pages 265-439, Springer
Berlin, 1986.

\bibitem{Wang2019}  X. Wang, {\em Blow-up solutions of the stochastic nonlocal heat equations},
Stoch. Dyn. {\bf19}  (2019) 1950014, 12 pp.

\bibitem{Yang2019}J. Yang, {\em Second critical exponent for a nonlinear nonlocal diffusion equation}, Appl. Math. Lett.
{\bf81}  (2018)  57-62.

\end{thebibliography}
 \end{document}